\documentclass[12pt,a4paper]{article}

\usepackage[T1]{fontenc}
\usepackage[latin1]{inputenc}
\usepackage[francais]{babel}
\usepackage{amsmath}
\usepackage{amssymb}
\numberwithin{equation}{section}

\newcommand{\K}{\mbox{$\mathbb{K}$}}
\newcommand{\C}{\mbox{$\mathbb{C}$}}
\newcommand{\N}{\mbox{$\mathbb{N}$}}

\newcommand{\Z}{\mbox{$\mathbb{Z}$}}
\newcommand{\zb}{\mbox{$\mathbb{Z}$}}
\newcommand{\Q}{\mbox{$\mathbb{Q}$}}
\newcommand{\F}{\mbox{$\mathbb{F}$}}

\def\ds{\displaystyle}
\def\inf{\infty}

\def\al{\alpha}
\def\ga{\gamma}

\def\de{\delta}

\def\la{\lambda}
\def\te{\theta}

\def\om{\omega}

\def\ca{\mathcal{A}}
\def\cb{\mathcal{B}}

\def\cg{\mathcal{G}}

\def\cs{\mathcal{S}}

\def\mb{\overline{m}}
\def\sark{S\'ark\"ozy }
\def\ni{\noindent}
\def \dv{\,|\,}
\def\lcm{\textrm{lcm}}

\begin{document}
\centerline{\bf On the parity of generalized partition functions
III}
\bigskip
\centerline{\bf by}
\bigskip

 \centerline {\bf F. Ben Sa\"id$^a$, J.-L.
Nicolas$^b$ and A. Zekraoui$^a$ \footnote {E-mail
addresses:Fethi.Bensaid@fsm.rnu.tn,jlnicola@in2p3.fr,ahlemzekraoui@yahoo.fr}
\footnote{Research partially supported by CNRS, by R\'egion
Rh\^one-Alpes, contract MIRA 2004 {\it Th\'eorie des nombres Lyon,
Saint-Etienne, Monastir} and by DGRST, Tunisia, UR 99/15-18.\\}}

\vskip 0.5 true cm
  \centerline {$^a$Univ. de Monastir, Facult\'e des Sciences, D\'epartement de Math\'ematiques,}
   \centerline { Avenue de l'environnement,
5000, Monastir, Tunisie.}
 \centerline { $^b$Institut Camile Jordan, UMR 5208,
Batiment Doyen Jean Braconnier,}
 \centerline { Univ. Claude Bernard (Lyon 1), 21
Avenue Claude Bernard, F-69622 Villeurbanne, France.}

\medskip
\vskip 1 true cm

\ni 
{\bf Abstract.} Improving on some results of J.-L. Nicolas
\cite {Ndeb}, the elements of the set ${\cal A}={\cal
A}(1+z+z^3+z^4+z^5)$, for which the partition function $p({\cal
A},n)$ (i.e. the number of partitions of $n$ with parts in ${\cal A}$)
is even for all $n\geq 6$ are determined. An asymptotic
estimate to the counting function of this set is also given.

\medskip
\ni

\medskip
\ni

\section {Introduction.}

Let $\mathbb{N}$ (resp. $\mathbb{N}_0$) be the set of positive
(resp. non-negative) integers. If ${\cal A}=\{a_1,a_2,...\}$ is a
subset of $\mathbb{N}$ and $n\in \mathbb{N}$ then $p({\cal A},n)$
is the number of partitions of $n$ with parts in $\cal A$, i.e.,
the number of solutions of the diophantine equation
\begin{equation}\label{11}
 a_1x_1+a_2x_2+\ldots=n,
\end{equation}
in non-negative integers $x_1,x_2,...$. As usual we set $p({\cal
A},0)=1$.\\ The counting function of the set ${\cal A}$ will be
denoted by $A(x)$, i.e.,
\begin{equation}\label{12}
A(x)=\mid\{n\leq x,~n\in {\cal A}\}\mid .
\end{equation}
Let $\F_2$ be the field with $2$ elements,
$P=1+\epsilon_1z^1+...+\epsilon_Nz^N\in \F_2[z],~N\geq 1.$
Although it is not difficult to prove (cf. \cite{NRS},
\cite{Ben2}) that there is a unique subset ${\cal A}={\cal A}(P)$
of $\mathbb{N}$ such that the generating function $F(z)$ satisfies
\begin{equation}\label{13}
F(z)=F_{\cal A}(z)=\prod_{a\in \ca }\frac {1}{1-z^a}=\sum_{n\geq
0}p({\cal A},n)z^n\equiv P(z) \pmod{2},
\end{equation}
the determination of the elements of such sets for general $P's$
seems to be hard. 

\medskip

Let the decomposition of $P$ into irreducible factors over $\F_2$
be
\begin{equation}\label{14}
P=P_1^{\alpha_1}P_2^{\alpha_2}...P_l^{\alpha_l}.
\end{equation}
We denote by $\beta_i=$ ord$(P_i)$, $1\leq i\leq l$, the order of $P_i$, 
that is the smallest positive integer $\beta_i$ such
that $P_i(z)$ divides $1+z^{\beta_i} $ in $\F_2[z]$. It is known
that $\beta_i $ is odd (cf. \cite{LN}). We set
\begin{equation}\label{15}
\beta = \lcm(\beta_1,\beta_2,...,\beta_l).
\end{equation}
Let ${\cal A}={\cal A}(P)$ satisfy (\ref {13}) and $\sigma ({\cal
A},n)$ be the sum of the divisors of $n$ belonging to ${\cal A}$,
i.e.,
\begin{equation}\label{16}
\sigma ({\cal A},n)=\sum_{d\mid n,\,d\in {\cal A}}d=\sum_{d\mid
n} d\chi({\cal A},d),
\end{equation}
where $\chi ({\cal A},.)$ is the characteristic function of the
set $\cal A$, i.e, $\chi ({\cal A},d)=1$ if $d\in \cal A$ and
$\chi ({\cal A},d)=0$ if $d\not \in \cal A$. It was proved in
\cite{Ben3} (see also \cite{Ben1},~\cite{Lah}) that for all $k\geq
0$, the sequence $(\sigma ({\cal A},2^kn)\bmod ~ 2^{k+1})_{n\geq
1}$ is periodic with period $\beta $ defined by (\ref {15}), in
other words,
\begin{equation}\label{17}
n_1\equiv n_2\pmod{\beta} \Rightarrow \forall k\geq 0,~\sigma
({\cal A},2^kn_1)\equiv \sigma ({\cal A},2^kn_2)\;(\bmod~2^{k+1}).
\end{equation}
Moreover, the proof of (\ref {17}) in \cite {Ben3} allows to
calculate $\sigma ({\cal A},2^kn)\bmod ~ 2^{k+1}$ and to deduce
the value of $ \chi({\cal A},n)$ where $n$ is any positive
integer. Indeed, let
\begin{equation}\label{18}
S_{\cal A}(m,k)=\chi({\cal A},m)+2\chi({\cal
A},2m)+\ldots+2^k\chi({\cal A},2^km).
\end{equation}
If $n$ writes $n=2^km$ with $k\ge 0$ and $m$ odd, (\ref{16})
implies
\begin{equation}\label{19}
\sigma ({\cal A}, n)=\sigma ({\cal A}, 2^km)=\sum_{d\dv m}
dS_{\cal A}(d,k),
\end{equation}
which, by M\"obius inversion formula, gives
\begin{equation}\label{110}
mS_{\cal A}(m,k)=\sum_{d\dv m} \mu(d) \sigma(\ca,\frac {
n}{d})=\sum_{d\dv \overline{m}} \mu(d) \sigma(\ca,\frac { n}{d}),
\end{equation}
where $\ds \;\overline{m}=\prod_{p\dv m}p\;$ denotes the radical
of $m$ with $\overline {1}=1$.\\In the above sums, $\frac { n}{d}$
is always a multiple of $2^k$, so that, from the values of
$\sigma(\ca,\frac { n}{d})$, by (\ref {110}), one can determine
the value of $S_{\cal A}(m,k)~\bmod ~2^{k+1}$ and by (\ref {18}),
the value of $\chi(\ca,2^im)$ for all $i$, $i\le k$.

\medskip

Let $\beta$ be an odd integer $\geq 3$ and $(\mathbb{Z}/{\beta
\mathbb{Z}})^*$ be the group of invertible elements modulo
$\beta$. We denote by $<2>$ the subgroup of $(\mathbb{Z}/{\beta
\mathbb{Z}})^*$ generated by $2$ and consider its action $\star $
on the set $\mathbb{Z}/{\beta \mathbb{Z}}$ given by $a\star x=ax$
for all $a\in ~<2>$ and $x\in \mathbb{Z}/{\beta \mathbb{Z}}$. The
quotient set will be denoted by  $(\mathbb{Z}/{\beta
\mathbb{Z}})/_{<2>}$ and the orbit of some $n$ in
$\mathbb{Z}/{\beta \mathbb{Z}}$ by $O(n)$. For $P\in
\mathbb{F}_2[z]$ with $P(0)=1$ and ord$(P)=\beta$, let ${\cal
A}={\cal A}(P)$ be the set obtained from (\ref {13}). Property
(\ref {17}) shows (after \cite{Bac1}) that if $ n_1$ and $ n_2$
are in the same orbit then
\begin{equation}\label{111}
\sigma ({\cal A},2^kn_1)\equiv \sigma ({\cal A},2^kn_2)~(\bmod~
2^{k+1}),~\forall k\geq 0.
\end{equation}
Consequently, for fixed $k$, the number of distinct values that
$(\sigma ({\cal A},2^kn)~\bmod ~ 2^{k+1})_{n\geq 1}$ can take is
at most equal to the number of orbits of $\mathbb{Z}/{\beta
\mathbb{Z}}$.

\medskip

Let $\varphi $ be the Euler function and $s$ be the order of $2$
modulo $\beta $, i.e., the smallest positive integer $s$ such that
$2^s\equiv 1~(\bmod ~ \beta )$. If $\beta =p $ is a prime number
then $(\mathbb{Z}/{p \mathbb{Z}})^*$ is cyclic and the number of
orbits of $\mathbb{Z}/{p \mathbb{Z}}$ is equal to $1+r$ with
$r=\frac {\varphi (p)}{s}=\frac {p-1}{s}$. In this case, we have
\begin{equation}\label{112}
(\mathbb{Z}/{p \mathbb{Z}})/_{<2>}=\{ O(g),~O(g^2),
...,~O(g^r)=O(1),~O(p)\},
\end{equation}
where $g$ is some generator of $(\mathbb{Z}/{p \mathbb{Z}})^*$.
For $r=2$, the sets ${\cal A}={\cal A}(P)$ were completely
determined by N. Baccar, F. Ben Sa\"id and J.-L. Nicolas
(\cite{Bac2}, \cite{Ben4}). Moreover, N. Baccar proved in
\cite{Bac3} that for all $r\geq 2$, the elements of ${\cal A}$ of
the form $2^km$, $k\geq 0$ and $m$ odd, are determined by the
$2$-adic development of some root of a polynomial with integer
coefficients. Unfortunately, his results are not explicit and do
not lead to any evaluation of the counting function of the set
${\cal A}$. When $r=6$, J.-L. Nicolas determined (cf. \cite{Ndeb})
the odd elements of ${\cal A}={\cal A}(1+z+z^3+z^4+z^5)$. His
results ( which will be stated in Section 2, Theorem 0) allowed to
deduce a lower bound for the counting function of ${\cal A}$. In
this paper, we will consider the case $p=31$ which satisfies
$r=6$. In $\F_2[z]$, we have
\begin{equation}\label{113}
\frac {1-z^{31}}{1-z}=P^{(1)}P^{(2)}...P^{(6)},
\end{equation}
with
$$P^{(1)}=1+z+z^3+z^4+z^5,~~P^{(2)}=1+z+z^2+z^4+z^5,P^{(3)}=1+z^2+z^3+z^4+z^5,$$
$$P^{(4)}=1+z+z^2+z^3+z^5,~~P^{(5)}=1+z^2+z^5,~~P^{(6)}=1+z^3+z^5.$$
In fact, there are other primes $p$ with $r=6$. For instance,
$p=223$ and $p=433$.

\medskip

In Section 2, for ${\cal A}={\cal A}(P^{(1)})$, we evaluate the
sum $S_{\cal A}(m,k)$ which will lead to results of Section 3
determining the elements of the set ${\cal A}$. Section $4$ will
be devoted to the determination of an asymptotic estimate to the
counting function $A(x)$ of $\cal A$. Although, in this paper,
the computations are only carried out for $P=P^{(1)}$, the results 
could probably be extended to any $P^{(i)}$, $1\leq i\leq 6$, and 
more generally, to any polynomial $P$ of order $p$ and such that $r=6$.

\medskip
{\bf Notation.} We write $a\bmod ~ b$ for the remainder of
the euclidean division of $a$ by $b$. The ceiling of the real number $x$
is denoted by
$$\lceil x\rceil ={\rm inf} \{n\in \mathbb{Z},~x\leq n\}.$$.

\vspace{-10mm}

\section {The sum $ S_{{\cal A}}(m,k),~{\cal A}={\cal A}(1+z+z^3+z^4+z^5).$}
From now on, we take ${\cal A}={\cal A}(P)$ with
\begin{equation}\label{21}
P=P^{(1)}=1+z+z^3+z^4+z^5.
\end{equation}
The order of $P$ is $\beta =31$. The smallest primitive root
modulo $31$ is $3$ that we shall use as a generator of
$(\mathbb{Z}/{31 \mathbb{Z}})^*$. The order of $2$ modulo $31$ is
$s=5$ so that
\begin{equation}\label{22}
(\mathbb{Z}/{31 \mathbb{Z}})/_{<2>}=\{ O(3),~O(3^2),
...,~O(3^6)=O(1),~O(31)\},
\end{equation}
with
\begin{equation}\label{23}
O(3^j)=\{2^k3^j,~0\leq k\leq 4\},~1\leq j\leq 6
\end{equation}
and
\begin{equation}\label{24}
O(31)=\{31n,~n\in \N \}.
\end{equation}  
For $k\geq 0$ and $0\leq j\leq 5$,  we define the integers $u_{k,j}$ by
\begin{equation}\label{25}
u_{k,j}=\sigma ({\cal A},2^k3^j)~\bmod ~ 2^{k+1}.
\end{equation}

{\bf The Graeffe transformation.} Let ${\K}$ be a field and
${\K}[[z]]$ be the ring of formal power series with coefficients
in ${\K}$. For an element
$$f(z)=a_0+a_1z+a_2z^2+\ldots+a_nz^n+\ldots$$ of this ring, the
product $$f(z)f(-z)=b_0+b_1z^2+b_2z^4+\ldots+b_nz^{2n}+\ldots$$ is
an even power series. We shall call $\cg(f)$ the series
\begin{equation}\label{Graeffe}
\cg(f)(z)=b_0+b_1z+b_2z^2+\ldots+b_nz^{n}+\ldots.
\end{equation}
It follows immediately from the above definition that for $f,~g\in {\K}[[z]]$,
\begin{equation}\label{g(fg)}
\cg(fg)=\cg(f)\cg(g).
\end{equation}
Moreover if $q$ is an odd integer and $f(z)=1-z^q$, we have
$\cg(f)=f$.
We shall use the following notation for the iterates
of $f$ by $\cg$:
\begin{equation}\label{iter}
f_{(0)}=f,\;\; f_{(1)}=\cg(f),\;\;\ldots,\;\;
f_{(k)}=\cg(f_{(k-1)})=\cg^{(k)}(f).
\end{equation}
More details about the Graeffe transformation are given in \cite{Ben3}.
By making the logarithmic derivative of formula (\ref {13}), we 
get (cf. \cite{NRS}):
\begin{equation}\label{derlog}
\sum_{n=1}^{\infty }\sigma ({\cal A},n)z^n=z\frac
{F'(z)}{F(z)}\equiv z\frac {P'(z)}{P(z)}~(\bmod~2),
\end{equation}
which, by Propositions 2 and 3 of \cite{Ben3}, leads to
\begin{equation}\label{sig1}
\sum_{n=1}^\infty \sigma ({\cal A},2^kn)z^n\equiv 
z\frac{P_{(k)}'(z)}{P_{(k)}(z)}= \frac{z}{1-z^{31}}
\left(P_{(k)}'(z)W_{(k)}(z)\right)~(\bmod~2^{k+1}),
\end{equation}
with $P_{(k)}'(z)=\frac{{\rm d}}{{\rm d}z}(P_{(k)}(z)$ and
\begin{equation}\label{W}
W(z)=(1-z)P^{(2)}(z)...P^{(6)}(z).
\end{equation}
Formula (\ref {sig1}) proves (\ref {111}) with $\beta =31$, and
the computation of the $k$-th iterates $P_{(k)}$ and $W_{(k)}$ by
the Graeffe transformation yields the value of $\sigma ({\cal
A},2^kn)~\bmod~2^{k+1}$. For instance, for $k=11$, we obtain:
$$u_{k,0}=1183,~~u_{k,1}=1598,~~u_{k,2}=1554,~~u_{k,3}=845,
~~u_{k,4}=264,~~u_{k,5}=701.$$
A divisor of $2^k3^j$ is either a divisor of $2^{k-1}3^j$ or 
a multiple of $2^k$. Therefore, from (\ref {25}) and (\ref{16}), 
$u_{k,j}\equiv u_{k-1,j} \pmod{2^{k}}$ holds
and the sequence $(u_{k,j})_{k\geq 0}$ defines a $2$-adic integer $U_j$
satisfying for all $k'$s:
\begin{equation}\label{Uicong}
U_j\equiv u_{k,j}~(\bmod ~2^{k+1}),~~0\leq j\leq 5.
\end{equation}
It has been proved in \cite{Bac3} that the $U_j'$s are the roots
of the polynomial $$R(y)=y^6-y^5+3y^4-11y^3+44y^2-36y+32.$$ Note
that $R(y)^5$ is the resultant in $z$ of
$\phi_{31}(z)=1+z+...+z^{30}$ and $y+z+z^2+z^4+z^8+z^{16}$.

\medskip

Let us set $$\theta =U_0=1+2+2^2+2^3+2^4+2^7+2^{10}+...$$ It turns
out that the Galois group of $R(y)$ is cyclic of order $6$ and
therefore the other roots $U_1,...,U_5$ of $R(y)$ are polynomials
in $\theta $. With Maple, by factorizing $R(y)$ on $\Q [\theta ]$
and using the values of $u_{11,j}$, we get 
$$U_0=\theta \equiv 1183 \pmod{2^{11}}$$  
$$U_1=\frac{1}{32}(3\theta^5+5\theta^3-36\theta^2+84\theta ) \equiv 1598
\pmod{2^{11}}$$ 
$$U_2=\frac{1}{32}(-3\theta^5-5\theta^3+20\theta^2-100\theta ) \equiv 1554
\pmod{2^{11}}$$ 
$$U_3=\frac{1}{32}(-\theta^5-7\theta^3+12\theta^2-44\theta +32)\equiv 845
\pmod{2^{11}}$$ 
$$U_4=\frac {1}{32}(-\theta^5+4\theta^4+\theta^3+24\theta^2-68\theta +96)
\equiv 264 \pmod{2^{11}}$$ 
\begin{equation}\label{Ui}
U_5=\frac {1}{16}(\theta^5-2\theta^4+3\theta^3-10\theta^2+48\theta
-48) \equiv 701 \pmod{2^{11}}.
\end{equation}
For convenience, if $j\in \Z$, we shall set
\begin{equation}\label{Uimod6}
U_j=U_{j \bmod ~6}.
\end{equation}

We define the completely additive function
$\ell:~\mathbb{Z}\setminus 31\mathbb{Z}\rightarrow
\mathbb{Z}/6\mathbb{Z}$ by
\begin{equation}\label{26}
\ell (n)=j~~\text ~~if~n\in O(3^j),
\end{equation}
so that $\ell (n_1n_2)\equiv \ell (n_1)+\ell (n_2)~(\bmod ~6)$. We
split the odd primes different from $31$ into six classes
according to the value of $\ell $. More precisely, for $0\leq
j\leq 5$,
\begin{equation}\label{27}
p\in {\cal P}_j\Longleftrightarrow \ell(p)=j \Longleftrightarrow
p\equiv 2^k3^j ~(\bmod ~31),~k=0,1,2,3,4.
\end{equation}
We take $L:~\mathbb{N}\setminus 31\mathbb{N}\longrightarrow
\mathbb{N}_0$ to be the completely additive function defined on
primes by
\begin{equation}\label{28}
L(p)=\ell (p).
\end{equation}
We define, for $0\leq j\leq 5$, the additive function
$\om_j:~\mathbb{N}\longrightarrow \mathbb{N}_0$ by
\begin{equation}\label{29}
\om_j(n)=\sum_{p\mid n,~p\in {\cal P}_j}1=\sum_{p\mid n,~\ell
(p)=j}1,
\end{equation}
and $\om(n)=\om_0(n)+...+\om_5(n)=\sum_{p\mid n} 1$. We remind that additive
functions vanish on $1$.

\medskip

From (\ref{25}), (\ref{23}), (\ref {111}) and (\ref{Uicong}), it
follows that if $n=2^km\in O(3^j)$ (so that $j=\ell(n)=\ell(m)$),
\begin{equation}\label{sigmaequivmod2k+1}
\sigma(\ca,n)=\sigma(\ca,2^km)\equiv U_{\ell(m)} \pmod{2^{k+1}}.
\end{equation}
We may consider the $2$-adic number
\begin{equation}\label{Sm}
S(m)=S_{\ca}(m)=\chi (\ca,m)+2\chi (\ca,2m)+...+2^k\chi(\ca,2^km)+...
\end{equation}
satisfying from (\ref{18}),
\begin{equation}\label{Smequivmod2k+1}
S(m)\equiv S_{\ca}(m,k) \pmod{2^{k+1}}.
\end{equation}
Then (\ref{110}) implies for $(m,31)=1$,
\begin{equation}\label{mSm}
mS(m)=\sum_{d\dv \overline{m}}\mu (d)U_{\ell (\frac {m}{d})}.
\end{equation}
If $31$ divides $m$, it was proved in \cite[(3.6)]{Bac1} that, for all
$k'$s,
\begin{equation}\label{Sm31}
\sigma (\ca,2^km)\equiv -5 \pmod{2^{k+1}}.
\end{equation}

{\bf Remark 1.} No element of $\ca $ has a prime factor in ${\cal
P}_0$. This general result has been proved in \cite{Bac1}, but we
recall the proof on our example: let us assume that 
$n=2^km\in {\cal A}$, where $m$ is an odd integer divisible by
some prime $p$ in ${\cal P}_0$, in other words $\om_0(m)\geq 1$.
(\ref {110}) gives
\begin{center}
\begin{eqnarray}
 mS_{\cal A}(m,k)&=&\sum_{d\dv m} \mu(d)
\sigma\left(\ca,\frac {n}{d}\right)=\sum_{d\dv \overline{m}}
\mu(d) \sigma\left(\ca,2^k\frac {m}{d}\right) \nonumber\\ & = &
\sum_{d\dv \frac {\overline{m}}{p}} \mu(d)
\sigma\left(\ca,2^k\frac {m}{d}\right)+\sum_{d\dv \frac
{{\overline{m}}}{p}} \mu(pd) \sigma\left(\ca,2^k\frac
{m}{pd}\right)\nonumber\\ & =&\sum_{d\dv \frac {\overline {
m}}{p}} \mu(d) \left( \sigma\left(\ca,2^k\frac {m}{d}\right)-
\sigma\left(\ca,2^k\frac { {m}}{pd}\right)\right).\nonumber
\end{eqnarray}
\end{center}
In the above sum, both $\frac {m}{d}$ and $\frac {m}{pd}$ are in
the same orbit, so that from (\ref {111}), $\sigma(\ca,2^k\frac
{m}{d})\equiv \sigma(\ca,2^k\frac {m}{pd})~(\bmod ~ 2^{k+1})$ and
therefore $mS_{\cal A}(m,k)\equiv 0~(\bmod ~ 2^{k+1})$. Since $m$
is odd and (cf. (\ref {18})) $0\leq S_{\cal A}(m,k)<2^{k+1}$ then
$S_{\cal A}(m,k)=0$, so that by (\ref{18}), $2^hm\not \in {\cal
A}$, for all $0\leq h\leq k$. 

\medskip
In \cite{Ndeb}, J.-L. Nicolas has
described the odd elements of ${\cal A}$. In fact, he obtained the
following:

\medskip

{\bf Theorem 0.} (\cite{Ndeb})) {\it (a) The odd elements of $\cal
A$ which are primes or powers of primes are of the form $p^\lambda
,~\lambda \geq 1$, satisfying one of the following four
conditions:
\begin{eqnarray}\label{Ndeb1}
p\in {\cal P}_1&and&\lambda \equiv 1,3,4,5~(\bmod
~6)\nonumber\\p\in {\cal P}_2 & and &\lambda \equiv 0,1~(\bmod~
3)\nonumber\\ p\in {\cal P}_4 & and & \lambda \equiv 0,1~(\bmod
~3)\nonumber\\ p\in {\cal P}_5 & and &\lambda \equiv
0,2,3,4~(\bmod ~6).\nonumber
\end{eqnarray}
(b) No odd element of $\cal A$ is a multiple of $31^2$. If $m$ is
odd, $m\not =1$, and not a multiple of $31$, then 
$$m\in {\cal A}~\text~if~and~only~if~~~~~~31m\in {\cal A}.$$ 
(c) An odd element
$n\in {\cal A}$ satisfies $\om_0(n)=0$ and $\om_3(n)=0$ or $1$; in
other words, $n$ is free of prime factor in ${\cal P}_0$ and has
at most one prime factor in ${\cal P}_3$.

(d) The odd elements
of ${\cal A}$ different from $1$, not divisible by $31$, which are
not primes or powers of primes are exactly the odd $n's$, $n\not
=1$, such that (where $\overline{n}=\prod_{p\mid n}p$):
\begin{enumerate}
  \item  $\om_0(n)=0$ and $\om_3(n)=0$ or $1$.
  \item  If $\om_3(n)=1$ then
$\ell (n)+\ell (\overline{n})\equiv 0~\text~or~1~(\bmod ~3).$
  \item If $\om_3(n)=0$ and $\om_1(n)+\ell (n)-\ell ( \overline{n})$ is
even then $$2\ell (n) -\ell ( \overline{n})\equiv 2\text
~or~3~or~4~or~5~~(\bmod ~6).$$
 \item If $\om_3(n)=0$ and $\om_1(n)+\ell
(n)-\ell ( \overline{n})$ is odd then 
$$2\ell (n) -\ell (\overline{n})\equiv 0\text ~or~4~~(\bmod ~6).$$
\end{enumerate}        }

{\bf Remark 2.} Point (b) of Theorem $0$ can be improved in the
following way: No element of $\ca $ is a multiple of $31^2$.
Indeed, from (\ref{110}), we have for $m$ odd, $k\geq 0$ and  $\tau
\geq 2$, 
$$31^\tau m S_{\ca }(31^\tau m,k)=\sum_{d\dv31^\tau m}\mu
(d)\sigma \left(\ca,2^k31^\tau \frac {m}{d}\right)=\sum_{d\dv
31\overline{m}}\mu (d)\sigma \left(\ca,2^k31^\tau \frac
{m}{d}\right)$$
$$=\sum_{d\dv \overline{m}}\mu (d) \left\{ \sigma
\left( \ca,2^k31^\tau \frac {m}{d}\right)-\sigma
\left(\ca,2^k31^{\tau -1}\frac {m}{d}\right) \right\}.$$ 
Since $31^\tau \frac {m}{d}$ and $31^{\tau -1}\frac {m}{d}$ are in the
same orbit $O(31)$ then (\ref {111}) and (\ref{Sm31}) give $\sigma
(\ca,2^k31^\tau \frac {m}{d})\equiv \sigma (\ca,2^k31^{\tau
-1}\frac {m}{d})\equiv -5 \pmod{2^{k+1}}$, so that we get $S_{\ca
}(31^\tau m,k)\equiv 0 \pmod{2^{k+1}}$. Hence, from (\ref {18}),
$S_{\ca }(31^\tau m,k)=0$ and for all $0\leq h\leq k$ and all $\tau \geq 2$,
$2^h31^\tau m$ does not belong to $\ca$.

\medskip 

In view of stating Theorem $1$ which will extend Theorem $0$, we
shall need some notation. The radical $\overline{m}$ of an odd
integer $m\ne 1 $, not divisible by $31$ and free of prime factors
belonging to ${\cal P}_0$ will be written
\begin{equation}\label{210}
\overline{m}=p_1\ldots p_{\om_1}p_{\om_1+1}\ldots p_{\om_1+\om_2}p_{\om_1+\om_2+1}
\ldots\ldots p_{\om_1+\om_2+\om_3+\om_4+1}\ldots p_{\om},
\end{equation}
where $\ell (p_i)=j$ for $\om_1+...+\om_{j-1}+1\leq i\leq
\om_1+...+\om_j$, $\om_j=\om_j(m)=\om_j(\overline{m})$ and
$\om=\om(m)=\om(\overline{m})\geq 1 $. We define the additive functions
from $\zb \setminus 31\zb$ into $\zb /12\zb $:
\begin{equation}\label{211}
\alpha =\alpha(m)= 2\om_5-2\om_1+\om_4-\om_2\bmod ~ 12,
\end{equation}
\begin{equation}\label{212}
a=a(m)=\om_5-\om_1+\om_2-\om_4\bmod ~ 12.
\end{equation}
Let $(v_i)_{i\in \mathbb{Z}}$ be the periodic sequence of period
$12$ defined by
\begin{equation}\label{vi} 
v_i=\left\{
\begin {array}{cl}
\frac {2}{\sqrt{3}}\cos (i\frac {\pi }{6})& \hbox{if $i$ is odd}\\
2\cos (i\frac {\pi }{6}) &\hbox{if $i$ is even.}
\end{array}\right.
\end{equation}
The values of $(v_i)_{i\in \mathbb{Z}}$ are given by:

\begin{center}
\begin{tabular}{|l|l|l|l|l|l|l|l|l|l|l|l|l|}
  \hline
 $i=$ & 0 & 1 & 2 & 3 & 4 & 5 & 6 & 7 & 8 & 9 & 10 & 11 \\
 $ v_i=$ & 2 & 1 & 1 & 0 & -1 & -1 & -2 & -1 & -1 & 0 & ~1 & ~1 \\
  \hline
\end{tabular}
\end{center}
Note that
\begin{equation}\label{vi6}
v_{i+6}=-v_i,
\end{equation}

\begin{equation}\label{vip2}
 v_i+v_{i+2}=\left\{
\begin {array}{cl}
v_{i+1}& \hbox{if $i$ is odd}\\
3v_{i+1} &\hbox{if $i$ is even,}
\end{array}\right.
\end{equation}
\begin{equation}\label{v2i}
v_{2i}\equiv -2^i~(\bmod ~3)
\end{equation}
and
\begin{equation}\label{vip6}
v_i\equiv v_{i+3}\equiv v_{2i} ~(\bmod ~2).
\end{equation}
From the $U_j$'s (cf. \eqref{Uicong} and \eqref{Ui}), we introduce 
the following $2$-adic integers:
\begin{equation}\label{E}
E_i=\sum_{j=0}^5v_{i+2j}U_j,~~i\in\Z,
\end{equation}
\begin{equation}\label{F}
F_i=\sum_{j=0}^5v_{i+4j}U_j,~~i\in\Z,
\end{equation}
\begin{equation}\label{G}
G=\sum_{j=0}^5(-1)^jU_j.
\end{equation}
From (\ref{vi6}), we have
\begin{equation}\label{Eiplus6}
E_{i+6}=-E_i,~~E_{i+12}=E_i,~~F_{i+6}=-F_i,~~F_{i+12}=F_i.
\end{equation}
From (\ref{vip2}), it follows that, if $i$ is odd,
\begin{equation}\label{Eiodd}
E_i+E_{i+2}=E_{i+1},~~F_i+F_{i+2}=F_{i+1},
\end{equation}
while, if $i$ is even,
\begin{equation}\label{Eieven}
E_i+E_{i+2}=3E_{i+1},~~F_i+F_{i+2}=3F_{i+1},
\end{equation}
The values of these numbers are given in the following array:

\begin{center}
\begin{tabular}{|l|l|l|}
  \hline
  $Z$ &  & $Z\bmod 2^{11}$ \\
  \hline
 $ E_0=$ & $\frac{1}{32} (11\theta^5-8\theta^4+29\theta^3
-124\theta^2+500\theta-256)$ & 1157 \\
  $E_1= $ & $\frac{1}{16} (3\theta^5-2\theta^4+9\theta^3
-26\theta^2+136\theta-64)$ & 1533 \\
   $E_2=$ & $3E_1-E_0$ & 1394 \\
   $E_3=$ & $2E_1-E_0$ & 1909 \\
  $E_4=$ & $3E_1-2E_0$ & 237 \\
  $E_5=$ & $E_1-E_0$ & 376 \\
    \hline
  $F_0=$ & $\frac{1}{32} (-3\theta^5-21\theta^3+36\theta^2-36\theta+64)$ & 1987 \\
  $F_1=$ & $\frac{1}{32} (-3\theta^5-4\theta^4-13\theta^3+24\theta^2-28\theta-64)$ & 166 \\
  $F_2=$ & $3F_1-F_0$ & 559 \\
  $F_3=$ & $2F_1-F_0$ & 393 \\
   $F_4=$ & $3F_1-2F_0$ & 620 \\
  $F_5=$ & $F_1-F_0$ & 227 \\
  \hline
  $G=$ & $\frac{1}{4} (-\theta^5+\theta^4-\theta^3+11\theta^2-34\theta+20)$ & 1905 \\
  \hline
\end{tabular}
\end{center}
\begin{center}
TABLE 1
\end{center}

{\bf Lemma 1.} {\it The polynomials $(U_j)_{0\leq j\leq 5}$ (cf. \eqref{Ui})
form a basis of $\Q [\theta ]$. The polynomials
$E_0,~E_1,~F_0,~F_1,~G,~U_0 $ form an other basis of $\Q
[\theta ]$. For all $i'$s, $E_i$ and $F_i$ are linear combinations
of respectively $E_0$ and $E_1$ and $F_0$ and $F_1$.}

\medskip

{\bf Proof.} With Maple, in the basis $1,\te,\ldots,\te^5$, we compute 
determinant $(U_0,...,U_5)=\frac {1}{1024}$. From \eqref{E}, \eqref{F} and
\eqref{G}, the determinant of $(E_0,E_1,$ $F_0,F_1,G,U_0 )$ in the basis 
$U_0,U_1,\ldots,U_5$ is equal to $12$. The last point
follows from (\ref {Eiodd}) and (\ref{Eieven}). $\Box$

\medskip

We have 

\medskip

{\bf Theorem 1.} {\it Let $m\ne 1$ be an odd integer not divisible
by $31$ with $\mb$ of the form (\ref {210}). Under the above notation and
the convention
\begin{equation}\label{0puiss}
 0^\om=\left\{
\begin {array}{cl}1& \hbox{if}~ \om=0
\\ 0&\hbox{if} ~\om>0,
\end{array}\right.
\end{equation}
we have:

\medskip
\ni
1) The $2$-adic integer $S(m)$ defined by (\ref {Sm}) satisfies 
$$mS(m)= 2^{\om_3-1}3^{\lceil \frac {
\om_2+\om_4}{2}-1\rceil}E_{\alpha  -2\ell (m)}+\frac
{0^{\om_3}}{2}3^{\lceil \frac { \om}{2}-1\rceil }
F_{a-4\ell (m)}$$
\begin{equation}\label{228}
 + \frac {0^{\om_2+\om_4}}{3}2^{\om-1}(-1)^{\ell (m)}G.
\end{equation}
2) The $2$-adic integer $S(31m)$ satisfies
\begin{equation}\label{231}
S(31m)= -31^{-1}S(m),
\end{equation}
where $31^{-1}$ is the inverse of $31$ in $\Z_2$. In particular, 
for all $k\in \{0,1,2,3,4\}$, we have 
$$2^km\in {\cal A}\quad \Longleftrightarrow \quad 
31\cdot 2^km\in {\cal A},$$ 
since the inverse of $31$ modulo $2^{k+1}$ is $-1$ for $k\leq 4$.}

\medskip

{\bf Proof of Theorem 1, 1).} From (\ref {mSm}), we have
\begin{equation}\label{232}
mS(m)=\sum_{d\dv \overline{m}} \mu(d) U_{\ell (\frac {m}{d})}
=\sum_{d\dv \overline{m}} \mu(d) U_{\ell(m)-\ell(d)}.
\end{equation}
Further, (\ref {232}) becomes
\begin{equation}\label{234}
mS(m)=\sum_{j=0}^5T(m,j)U_{\ell (m)-j}=\sum_{j=0}^5T(m,\ell
(m)-j)U_{j},
\end{equation}
with
\begin{equation}\label{235}
T(m,j)=T(\overline{m},j)=\sum_{d\mid \overline{m},~\ell (d)\equiv
j~(\bmod ~ 6)}\mu (d).
\end{equation}
Therefore (\ref {228}) will follow from (\ref {234}) and from the
following lemma:

\medskip

{\bf Lemma 2.} {\it The integer $T(m,j)$ defined in (\ref {235})
with the convention (\ref {0puiss}) and the definitions (\ref{29})
and (\ref {210})-(\ref {vi}), for $m\ne 1$, is equal to
$$T(m,j)=2^{\om_3-1}3^{\lceil \frac {\om_2+\om_4}{2}-1\rceil }v_{\alpha
-2j}+\frac {0^{\om_3}}{2}3^{\lceil \frac
{\om}{2}-1\rceil }v_{a -4j}$$
\begin{equation}\label{235bis}
+0^{\om_2+\om_4}\frac {(-1)^j}{3}2^{\om-1}.
\end{equation}
}

{\bf Proof.} Let us introduce the polynomial
\begin{equation}\label{236}
f(X)=(1-X)^{\om_1}(1-X^2)^{\om_2}...(1-X^5)^{\om_5}=\sum_{\nu \geq
0}f_\nu X^\nu .
\end{equation}
If the five signs were plus instead of minus, $f(X)$ would be the
generating function of the partitions in at most $\om_1$ parts equal
to $1$, ..., at most $\om_5$ parts equal to $5$. More generally, the
polynomial $$\widetilde{f}(X)=\prod
_{i=1}^\om(1+a_iX^{b_i})=\sum_{\nu \geq 0}\widetilde{f_\nu}X^\nu$$
is the generating function of
$$\widetilde{f_\nu}=\sum_{\epsilon_1, ..., \epsilon_\om\in
\{0,1\},~\sum_{i=1}^\om \epsilon_ib_i=\nu
}~~\prod_{i=1}^\om a_i^{\epsilon_i}.$$ 
To the vector
$\underline{\epsilon }=(\epsilon_1, ..., \epsilon_\om)\in \F_2^\om$,
we associate 
$$d=\prod _{i=1}^\om p_i^{\epsilon_i},~~
\mu(d)=\prod_{i=1}^\om (-1)^{\epsilon_i},~~
L(d)=\sum_{i=1}^\om\epsilon_i\ell (p_i)$$ 
where $L$ is the arithmetic function defined by \eqref{28} and we get
\begin{equation}\label{237}
f_\nu =\sum_{d\mid \overline{m},~L(d)=\nu}\mu (d),
\end{equation}

Consequently, by setting $\xi =\exp(\frac {i\pi }{3})$, (\ref
{235}), (\ref{236}) and (\ref {237}) give $$ T(m,j)=\sum_{\nu,~\nu
\equiv j~(\bmod ~ 6)}~~\sum_{d\mid \overline{m},~L(d)=\nu}\mu
(d)$$$$=\sum_{\nu \equiv j~(\bmod ~6)}f_\nu =\frac
{1}{6}\sum_{i=0}^5\xi^{-ij}f(\xi^i)=\frac
{1}{6}\sum_{i=1}^5\xi^{-ij}f(\xi^i)$$
\begin{equation}\label{238}
=\frac
{1}{6}\sum_{i=1}^5\xi^{-ij}(1-\xi^i)^{\om_1}(1-\xi^{2i})^{\om_2}(1-\xi^{3i})^{\om_3}(1-\xi^{4i})^{\om_4}(1-\xi^{5i})^{\om_5}.
\end{equation}
By observing that $$1-\xi =\xi^5,~1-\xi^2=\varrho
=\sqrt{3}(\cos\frac {\pi}{6}-i\sin \frac {\pi
}{6}),~1-\xi^3=2,~1-\xi^4 =\overline{\varrho},~1-\xi^6=0,$$ the
sum of the terms in $i=1$ and $i=5$ in (\ref {238}), which are
conjugate, is equal to
\begin{equation}\label{239}
\frac {2}{6} {\cal
R}(\xi^{-j}\xi^{5\om_1}\varrho^{\om_2}2^{\om_3}\overline{\varrho}^{\om_4}\xi^{\om_5})=\frac
{2^{\om_3}}{3}\sqrt{3}^{\om_2+\om_4}\cos \frac {\pi
}{6}(2\om_5-2\om_1+\om_4-\om_2-2j).
\end{equation}
Now, the contribution of the terms in $i=2$ and $i=4$ is

\begin{eqnarray}\label{240}
\frac {2}{6} {\cal R}(\xi^{-2j}\varrho^{\om_1}\overline{\varrho}^{\om_2}
0^{\om_3}\varrho^{\om_4}\overline{\varrho}^{\om_5})\hspace{-3mm}
&=&\hspace{-3mm} 0^{\om_3}\frac {\sqrt{3}^{\om_1+\om_2+\om_4+\om_5 }}{3}
\cos \frac {\pi}{6}(\om_2+\om_5-\om_1-\om_4-4j)\notag\\
&=& \hspace{-3mm}0^{\om_3}\frac {\sqrt{3}^{\om}}{3}
\cos \frac {\pi}{6}(\om_2+\om_5-\om_1-\om_4-4j)
\end{eqnarray}

Finally, the term corresponding to $i=3$ in (\ref {238}) is equal
to
\begin{equation}\label{241}
\frac {1}{6}
(-1)^j2^{\om_1}0^{\om_2}2^{\om_3}0^{\om_4}2^{\om_5}
=0^{\om_2+\om_4}\frac{(-1)^j}{6}2^{\om_1+\om_3+\om_5}
=0^{\om_2+\om_4}\frac{(-1)^j}{6}2^{\om}.
\end{equation}
Consequently, by using our notation (\ref {210})-(\ref {212}),
(\ref {238}) becomes 
$$ T(m,j)=\frac
{2^{\om_3}}{3}\sqrt{3}^{\om_2+\om_4}\cos \frac {\pi }{6}(\alpha -2j)
+0^{\om_3}\frac {\sqrt{3}^{\om}}{3}\cos \frac {\pi
}{6}(a-4j)$$
\begin{equation}\label{242}
+0^{\om_2+\om_4}\frac {(-1)^j}{6}2^{\om}.
\end{equation}
Observing that $\alpha -2j$ has the same parity than $\om_2+\om_4$ and
similarly for $a-4j$ and $\om$ (when $\om_0=\om_3=0$), via
(\ref{vi}), we get (\ref {235bis}).

\medskip

{\bf Proof of Theorem 1, 2).} For all $k\geq 0$, from (\ref{110}),
we have
\begin{eqnarray}\label{254}
31mS_{\cal A}(31m,k)&=&\sum_{d\dv 31m}\mu (d)\sigma (\ca ,31\cdot
2^k\frac {m}{d})=\sum_{d\dv 31\overline{m}}\mu (d)\sigma (\ca
,31\cdot 2^k\frac {m}{d})\nonumber\\ & = & \sum_{d\dv
\overline{m}}\mu (d)\sigma (\ca ,31\cdot 2^k\frac
{m}{d})-\sum_{d\dv \overline{m}}\mu (d)\sigma (\ca ,2^k\frac
{m}{d})\nonumber\\ &=&\sum_{d\dv \overline{m}}\mu (d)\sigma (\ca
,31\cdot 2^k\frac {m}{d})-mS_{\cal A}(m,k).
\end{eqnarray}
Since for all $d$ dividing $\overline{m}$, $31\cdot 2^k\frac
{m}{d}\in O(31)$ then, from (\ref {Sm31}), $\sigma (\ca ,31\cdot
2^k\frac {m}{d})\equiv \sigma(\ca,31\cdot2^k)\equiv -5\pmod
{2^{k+1}}$, so that (\ref {254}) gives 
\begin{equation}\label{255}
31mS_{\cal A}(31m,k)+mS_{\cal A}(m,k)\equiv -5\sum_{d\dv \overline{m}}\mu (d)
\pmod{2^{k+1}}.
\end{equation}
Since $\mb \neq 1$, $31mS_{\cal A}(31m,k)+mS_{\cal A}(m,k)\equiv 0 
\pmod{2^{k+1}}$. Recalling that $m$ is odd, by using
(\ref{Sm}), (\ref{Smequivmod2k+1}) and their similar for $S(31m)$,
we obtain the desired result. $\Box$

\section{Elements of the set $ {\cal A}={\cal A}(1+z+z^3+z^4+z^5).$}
 
In this section, we will determine the
elements of the set $\cal A$ of the form $n=2^k31^\tau m$, where
$\mb\ne 1$ satisfies (\ref {210}) and $\tau \in \{ 0,1\}$, since from
Remark 2, $2^k31^\tau m\notin \ca$ for all $\tau \geq 2$ . The
elements of the set
${\cal A}(1+z+z^3+z^4+z^5)$ of the form
$31^\tau 2^k$, $\tau =0$ or $1$, were shown in \cite {Bac3} to be
solutions of $2$-adic equations. More precisely, the following was
proved in that paper.

\medskip
\ni
1) The elements of the set ${\cal
A}(1+z+z^3+z^4+z^5)$ of the form $2^k$, $k\ge 0$, are given
by the $2$-adic solution 
$$\sum_{k\ge 0} \chi(\ca,2^k)\,2^k=S(1)
=U_0=1+2+2^2+2^3+2^4+2^7+2^{10}+2^{11}+...$$ 
of the equation 
$$y^6-y^5+3y^4-11y^3+44y^2-36y+32=0.$$ 
Note that $S(1)=U_0$ follows from \eqref{mSm}.

\medskip
\ni
2) The elements of the set ${\cal A}(1+z+z^3+z^4+z^5)$
of the form $31\cdot 2^k$, $k\ge 0$, are given by the
solution 
$$\sum_{k\ge 0} \chi(\ca,31\cdot 2^k)\,2^k=S(31)=y=2^2+2^5+2^{11}+...$$ 
of the equation
$$31^5y^6+31^5y^5+13\cdot 31^4y^4+91\cdot 31^3y^3+364\cdot
31^2y^2+796\cdot 31y+752=0,$$ 
since, from \eqref{255} with $m=1$, we have $31S(31)=-5-U_0$, so that
$$S(31)=\frac{5+U_0}{1-32}=(1+4+U_0)(1+2^5+2^{10}+...)=2^2+2^5+2^{11}+...$$

{\bf Theorem 2.} {\it Let $m\ne 1$ be an odd integer not divisible by
any prime $p\in {\cal P}_0$ (cf. (\ref{27})) neither by $31^2$.
Then the sum $S(m)$ defined by (\ref{Sm}) does not vanish. So we
may introduce the $2$-adic valuation of $S(m)$:
\begin{equation}\label{gammam}
\gamma =\gamma (m)=v_2(S(m)).
\end{equation}
Then, if $31$ does not divide $m$, we have
\begin{equation}\label{gamma31}
\gamma (31m)=\gamma (m).
\end{equation}
Let us assume now that $m$ is coprime with $31$. We shall use the
quantities $\om_i=\om_i(m)$ defined by (\ref{29}), $\ell (m)$, $\alpha
=\alpha (m)$, $a=a(m)$ defined by (\ref{26}), (\ref{211}) and
(\ref{212}),
\begin{equation}\label{alphaprime}
\alpha' =\alpha'(m)=\alpha -2\ell (m) \bmod
~12=2\om_5-2\om_1+\om_4-\om_2-2\ell (m) \bmod ~12,
\end{equation}
\begin{equation}\label{aprime}
a'=a'(m)=a-4\ell (m) \bmod ~12=\om_5-\om_1+\om_2-\om_4-4\ell (m) \bmod
~12,
\end{equation}
$$ t=t(m)=\left\lceil \frac {\om_1+\om_5+\om_2+\om_4}{2}-1\right\rceil -
\left\lceil \frac{\om_2+\om_4}{2}-1\right\rceil $$

\begin{equation}\label{tm}
 =\left\{
\begin {array}{crll}\lceil \frac {\om_1+\om_5}{2}\rceil & \hbox {if}&~ \om_1+\om_5\equiv ~\om_2+\om_4 \equiv 1~(\bmod
~2)
\\ \lceil \frac {\om_1+\om_5}{2}-1\rceil  &\hbox {if} & {not.}
\end{array}\right.
\end{equation}
We have:

{\bf (i)} if $\om_3\ne 0$ and $\om_2+\om_4\ne 0$, the value of $\gamma
=\gamma (m)$ is given by

$$ \gamma = \left\{
\begin {array}{crll}\om_3-1& \hbox{if}&~ \alpha'\equiv 0,1,3,4 ~(\bmod
~6)
\\ \om_3&\hbox{if}& ~\alpha'\equiv 2 ~(\bmod ~6)
\\ \om_3+2&\hbox{if} &~\alpha'\equiv 5 ~(\bmod ~6).

\end{array}\right.
$$

{\bf (ii)} If $\om_2+\om_4=0$ and $\om_3\geq 1$, we set 
$\al''=\al'+6\ell(m) \bmod 12$ and $\de(i)=v_2(E_i+2^{v_2(E_i)}G)$ and we have
\begin{eqnarray*} 
&\text{if}& \om_1+\om_5 < v_2(E_{\al''}), \quad \text{ then }
\quad \ga=\om_3-1+\om_1+\om_5,\\
&\text{if}& \om_1+\om_5 = v_2(E_{\al''}), \quad \text{ then }
\quad \ga=\om_3-1+\de(\al''),\\
&\text{if}& \om_1+\om_5 > v_2(E_{\al''}), \quad \text{ then }
\quad \ga=\om_3-1+v_2(E_{\al''}).
\end{eqnarray*}

{\bf (iii)} If $\om_3=0$ and $\om_2+\om_4\ne 0$, we have 
$$\gamma
=-1+v_2(E_{\alpha'}+3^tF_{a'}).$$

{\bf (iv)} If $\om_3=\om_2=\om_4=0$ and $\om_1+\om_5\ne 0$, we have 
$$\gamma=-1+v_2(E_{\alpha'}+3^tF_{a'}+2^{\om_1+\om_5}(-1)^{\ell(m)}G).$$ }

{\bf Proof.} We shall prove that $S(m)\ne 0$ in each of the four
cases above. Assuming $S(m)\ne 0$, it follows from Theorem 1, 2)
that $S(31m)\ne 0$ and that $\gamma (31m)=\gamma (m)$, which sets
(\ref{gamma31}).

\medskip

{\bf Proof of Theorem 2 (i).} In this case, formula
(\ref{228}) reduces to $$mS(m)=2^{\om_3-1}3^{\lceil \frac
{\om_2+\om_4}{2}-1\rceil }E_{\alpha'}.$$ Since $E_{\alpha'}\ne 0$,
$S(m)$ does not vanish; we have $$\gamma
=v_2(S(m))=\om_3-1+v_2(E_{\alpha'})$$ and the result follows from
the values of $E_{\alpha'}$ modulo $2^{11}$ given in Table $1$.

\medskip

{\bf Proof of Theorem 2 (ii).} If $\om_2+\om_4=0$ and $\om_3\ne 0$,
formula (\ref{228}) becomes (since, cf. \eqref{Eiplus6}, $E_{i+6}=-E_i$ holds) 
$$mS(m)\!
=\!\frac {2^{\om_3-1}}{3}\left(E_{\al'}+2^{\om_1+\om_5}(-1)^{\ell(m)}G\right)
\!=\!(-1)^{\ell(m)}\frac {2^{\om_3-1}}{3}\left(E_{\al''}+2^{\om_1+\om_5}G\right).$$ 
As displaid in Table 1, $E_i$ is a linear combination of $E_0$ and $E_1$ so 
that, from Lemma 1, $S(m)$ does not vanish and $\gamma
=\om_3-1+v_2\left(E_{\al''}+2^{\om_1+\om_5}G\right)$, whence the result.
The values of $v_2(E_i)$ and $\de(i)$ calculated from Table 1 are given below.
$$
\begin{array}{|c|c|c|c|c|c|c|c|c|c|c|c|c|}
\hline
i&0&1&2&3&4&5&6&7&8&9&10&11\\
\hline
v_2(E_i)&0&0&1&0&0&3&0&0&1&0&0&3\\
\hline
\de(i)&1&1&2&1&1&8&2&2&4&2&2&4\\
\hline
\end{array}
$$

\medskip

{\bf Proof of Theorem 2 (iii).} If $\om_3= 0$ and $\om_2+\om_4\ne 0$ it
follows, from (\ref{228}) and the definition of $t$ above, that
$$mS(m)=\frac {1}{2}~3^{\lceil \frac {\om_2+\om_4}{2}-1\rceil
}(E_{\alpha'}+3^tF_{a'}).$$ 
But $E_i$ and $F_i$ are non-zero linear combinations of, respectively,
$E_0$ and $E_1$ and $F_0$ and $F_1$; by Lemma 1, $E_{\al'}+3^tF_{a'}$ does 
not vanish and $\gamma =-1+v_2(E_{\alpha'}+3^tF_{a'})$.

\medskip

{\bf Proof of Theorem 2 (iv).} If $\om_3=\om_2=\om_4=0$ and $m\ne 1$,
formula (\ref{228}) gives 
$$mS(m)=\frac {1}{6}\left(E_{\alpha'}+
3^tF_{a'}+2^{\om_1+\om_5}(-1)^{\ell(m)}G\right).$$ 
From Lemma $1$, we obtain 
$E_{\alpha'}+3^tF_{a'}+2^{\om_1+\om_5}(-1)^{\ell(m)}G\not=0$, which
implies $S(m)\ne 0$ and $\gamma =-1+v_2\left(
E_{\alpha'}+3^tF_{a'}+2^{\om_1+\om_5}(-1)^{\ell(m)}G\right)$. $\Box$

{\bf Theorem 3.} {\it Let $m$ be an odd integer satisfying $m\ne
1$, $(m,31)=1$, and with $\mb$ of the form (\ref{210}). Let $\gamma =\gamma
(m)$ as defined in Theorem $2$ and $Z(m)$ be the odd part of the
right hand-side of (\ref{228}), so that
\begin{equation}\label{mSmfinal}
mS(m)=2^{\gamma (m)}Z(m).
\end{equation}
\indent
{\bf (i)} If $k<\gamma$, then $2^km\notin \ca$ and $2^k31m\notin
\ca$.

{\bf (ii)} If $k=\gamma$, then $2^km\in \ca$ and $2^k31m\in
\ca$.

{\bf (iii)} If $k=\gamma +r$, $r\geq 1$, then we set 
$\cs_r=\{ 2^r+1,2^r+3,...,2^{r+1}-1\}$ and we have
$$2^{\gamma+r} m\in{\ca }\quad\Longleftrightarrow\quad \exists ~ l\in
\cs_r,~m\equiv l^{-1}Z(m) \pmod{2^{r+1}},$$ 
$$2^{\gamma +r} 31m\in{\ca }\quad\Longleftrightarrow \quad\exists ~
l\in \cs_r,~m\equiv -(31l)^{-1}Z(m) \pmod{2^{r+1}}.$$    } 

{\bf Proof of Theorem 3, (i).} We remind that $m$
is odd and (cf. \ref{Smequivmod2k+1}) $S(m)\equiv S_{\ca
}(m,k)(\bmod ~2^{k+1})$. It is obvious from (\ref {mSmfinal}) that
if $\gamma >k$ then $ S_{\ca }(m,k)\equiv 0(\bmod ~2^{k+1})$. So
that from (\ref {18}), $ S_{\ca }(m,k)=0$ and $2^hm\not\in \ca $,
for all $h,~0\leq h\leq k$. To prove that $2^k31m\notin \ca$, it
suffices to use this last result and (\ref{231}) modulo
$2^{k+1}$.

\medskip

{\bf Proof of Theorem 3, (ii).} If $\gamma =k$ then the same
arguments as above show that $$ mS_{\ca }(m,k)\equiv 2^kZ(m)
(\bmod ~2^{k+1}).$$ So that, by using Theorem 3, (i) and
(\ref{18}), we obtain $$2^km\chi ({\ca },2^km)\equiv 2^kZ(m)
(\bmod ~2^{k+1}).$$ Since both $m$ and $Z(m)$ are odd, we get
$\chi ({\ca },2^km)\equiv 1 (\bmod ~2)$, which shows that $2^km\in
{\ca }$. Once again, to prove that $2^k31m\in \ca$, it suffices to
use this last result and (\ref{231}) modulo $2^{k+1}$.

\medskip

{\bf Proof of Theorem 3, (iii).} Let us set $k=\gamma +r$, $r\geq
1$. (\ref{mSmfinal}) and (\ref{Smequivmod2k+1}) give
\begin{equation}\label{Theorem3,(iii)}
 mS_{\ca }(m,k)\equiv 2^\gamma Z(m)
(\bmod ~2^{\gamma +r+1}).
\end{equation}
So that, by using Theorem 3, (i) and (ii), we get
$$ m(2^{\gamma }+2^{\gamma +1}\chi({\cal
A},2^{\gamma +1}m)+\ldots+2^{\gamma +r}\chi({\cal A},2^{\gamma +r}
m)) \equiv 2^{\gamma }Z(m)(\bmod ~2^{\gamma +r+1}),$$ 
which reduces to 
$$ m(1+2\chi({\cal
A},2^{\gamma +1}m)+\ldots+2^{r}\chi({\cal A},2^{\gamma +r} m))
\equiv Z(m)(\bmod ~2^{r+1}).$$ 
By observing that $2^{\gamma
+r}m\in \ca $ if and only if $l=1+2\chi({\cal A},2^{\gamma
+1}m)+\ldots+2^{r}\chi({\cal A},2^{\gamma +r} m)$ is an odd
integer in $\cs_r$, we obtain $$2^{\gamma
+r} m\in{\ca }\quad\Longleftrightarrow\quad m\equiv l^{-1}Z(m) \pmod {2^{r+1}},
~l\in \cs_r.$$ 
To prove the similar result for
$2^{\gamma +r} 31m$, one uses the same method and (\ref {231})
modulo $2^{k+1}$. $\Box$

\section{The counting function.}
In Theorem 4 below, we will determine an
asymptotic estimate to the counting function  $A(x)$ (cf. (\ref
{12})) of the set ${\cal A}={\cal A}(1+z+z^3+z^4+z^5)$. The
following lemmas will be needed.

\medskip

{\bf Lemma 3.} {\it Let $K$ be any positive integer and $x\geq 1$
be any real number. We have $$ \mid \{ n\leq x:~\gcd(n,K)=1\}\mid
\leq 7\frac {\varphi (K)}{K}x,$$ where $\varphi $ is the Euler
function.}

\medskip

{\bf Proof.} This is a classical result from sieve theory: see
Theorems $3-5$ of \cite {Hal}. $\Box$

\medskip

{\bf Lemma 4.} {\it (Mertens's formula) Let $\theta $ and $\eta $
be two positive coprime integers. There exists an absolute
constant $C_1$ such that, for all $x>1$, $$\pi (x;\theta ,\eta
)=\prod_{p\leq x,~p\equiv \theta (\bmod ~ \eta )}(1-\frac
{1}{p})\leq \frac {C_1}{(\log x)^{\frac {1}{\varphi (\eta )}}}.$$
}

{\bf Proof.} For $\theta $ and $\eta $ fixed, Mertens's formula
follows from the Prime Number Theorem in arithmetic progressions.
It is proved in \cite{BSLN} that the constant $C_1$ is absolute. $\Box$

\medskip

{\bf Lemma 5.} {\it For $i\in \{2,3,4\}$, let
$$K_i=K_i(x)=\prod_{p\leq x,~\ell (p)\in \{0,i\}}p=\prod_{p\leq
x,~ p\in {\cal P}_0\cup {\cal P}i}p,$$ where $\ell $, ${\cal P}_0$
and ${\cal P}_i$ are defined by (\ref {26})-(\ref{27}). Then for
$x$ large enough, $$\mid \{n: 1\leq n\leq x,~\gcd(n,K_i)=1\}\mid
={\cal O}\left(\frac {x}{(\log x)^{\frac {1}{3}}}\right).$$ }

{\bf Proof.} By Lemma $3$ and (\ref {27}), we have 
$$\mid \{
n:n\leq x,~\gcd(n,K_i)=1\}\mid \leq 7x\frac {\varphi
(K_i)}{K_i}$$
$$=7x\prod_{0\leq j\leq 4,\;\tau \in \{0,i\}}
\quad\prod_{p\leq x,\;p\equiv 2^j3^\tau (\bmod ~ 31)}(1-\frac
{1}{p}).$$ 
So that by Lemma $4$, for all $i\in \{2,3,4\}$ and $x$
large enough, 
$$ \mid \{ n:n\leq x,~\gcd(n,K_i)=1\}\mid \leq \frac
{7C_1^{10}x}{(\log x)^{\frac {10}{\varphi (31)}}}={\cal
O}\left(\frac {x}{(\log x)^{\frac {1}{3}}}\right). \;\Box$$

{\bf Lemma 6.} {\it Let $r,u\in \N_0$, $\ell$ and $\alpha' $ be the
functions defined by (\ref {26}) and (\ref {alphaprime}), $\om_j$ be
the additive function given by (\ref {29}). We take $\xi $ to be a
Dirichlet character modulo $2^{r+1}$ with $\xi_0$ as principal
character and we let $\varrho $ be the completely multiplicative
function defined on primes $p$ by
\begin{equation}\label{41}
\varrho (p)=\left\{
\begin {array}{crll}0& \hbox{if}~ \ell (p)=0 ~\hbox{or}~ p=31
\\ 1 &
\hbox{otherwise.}
\end{array}\right.
\end{equation}
If $y$ and $z$ are respectively some $2^u$-th and $12$-th roots of
unity in $\C$, and if $x$ is a real number $> 1$, we set
\begin{equation}\label{42}
S_{y,z,\xi }(x)=\sum_{~2^{\om_3(n)}n\leq x}\varrho (n)\xi
(n)y^{\om_2(n)+\om_4(n)} z^{\alpha'(n)}.
\end{equation}
Then, when $x$ tends to infinity, we have 

\ni
$\bullet $ If $\xi \ne \xi_0$,
\begin{equation}\label{43}
S_{y,z,\xi }(x)={\cal O}\left(x\frac {\log \log x}{(\log
x)^2}\right).
\end{equation}
$\bullet $ If $\xi = \xi_0$,
 \begin{equation}\label{45}
 S_{y,z,\xi_0 }(x)=\frac {x}{(\log x)^{1-f_{y,z}(1)}}
 \left( \frac {H_{y,z,\xi_0}(1)C_{y,z}}{\Gamma \left(f_{y,z}(1)\right)}+
  {\cal O} \left(\frac {\log \log x}{\log x}\right)\right),
\end{equation}
where $\Gamma $ is the Euler gamma function,
 \begin{equation}\label{46}
f_{y,z}(s)=\frac {1}{\varphi (31)}\sum_{1\leq j\leq
5}~\sum_{p,~\ell (p)=j}g_{j,y,z}(s),
\end{equation}
\begin{equation}\label{49}
g_{1,y,z}(s)=z^8,~g_{2,y,z}(s)= yz^7,~g_{3,y,z}(s)=\frac
{z^6}{2^s},~g_{4,y,z}(s)= yz^5,~g_{5,y,z}(s)=z^4,
\end{equation}
\begin{equation}\label{47}
H_{y,z,\xi }(s)=\prod_{1\leq j\leq 5}~\prod_{p,~\ell
(p)=j}\left(1+\frac {g_{j,y,z}(s)\xi
(p)}{p^s-z^{-2j}\xi (p)}\right)\left(1-\frac {\xi
(p)}{p^s}\right)^{g_{j,y,z}(s)},
\end{equation}
\begin{equation}\label{410}
C_{y,z}=\prod_{1\leq j\leq 5}\left\{ \prod_{p,~\ell (p)=j}(1-\frac
{1}{p})^{-g_{j,y,z}(1)}\prod_p(1-\frac {1}{p})^{\frac
{g_{j,y,z}(1)}{30}}\right\}.
\end{equation}         }

{\bf Proof.} The evaluation of such sums is based, as we know, on
the Selberg-Delange method. In \cite{Ben5}, one finds an
application towards direct results on such problems. In our case,
to apply Theorem 1 of that paper, one should start with expanding,
for complex number $s$ with ${\cal R} s>1$, the Dirichlet series
$$F_{y,z,\xi }(s)=\sum_{n\geq 1}\frac {\varrho (n)\xi
(n)y^{\om_2(n)+\om_4(n)} z^{\alpha' (n)}}{(2^{\om_3(n)}n)^s}$$ in an
Euler product given by
\begin{eqnarray}\label{Eulerprod}
F_{y,z,\xi }(s) & = & \prod_{1\leq j\leq 5}~\prod_{p,~\ell
(p)=j}\left( 1+\sum_{m=1}^\infty \frac { \xi (p^m)
y^{\om_2(p^m)+\om_4(p^m)}z^{\alpha'(p^m)}}{\left(2^{\om_3(p^m)}p^m\right)^s}\right)\nonumber\\
&=&\prod_{1\leq j\leq 5}~\prod_{p,~\ell (p)=j}\left( 1+\frac
{g_{j,y,z}(s)\xi (p)}{p^s-z^{-2j}\xi
(p)}\right),\nonumber
\end{eqnarray}
which can be written
$$F_{y,z,\xi }(s)=H_{y,z,\xi }(s)\prod_{1\leq j\leq
5}~\prod_{p,~\ell (p)=j}\left(1-\frac {\xi
(p)}{p^s}\right)^{-g_{j,y,z}(s)},$$ 
where $g_{j,y,z}(s)$ and $H_{y,z,\xi }(s)$ are defined by (\ref{49})
and (\ref{47}). To complete the proof of Lemma $6$, one
has to show that $H_{y,z,\xi }(s)$ is holomorphic for ${\cal R}
s>\frac {1}{2}$ and, for $y$ and $z$ fixed, that $H_{y,z,\xi }(s)$
is bounded for ${\cal R} s\geq \sigma_0 >\frac {1}{2}$, which can
be done by adapting the method given in \cite {Ben5} (Preuve du
Theor\`eme 2, p. 235). $\Box$

\medskip

{\bf Lemma 7.} {\it We keep the above notation and we let ${\cal
G}$ be the set of integers of the form $n=2^{\om_3(m)}m$ with the
following conditions:
\begin{itemize}
\item $m$ odd and $\gcd(m,31)=1$,
\item $m=m_1m_2m_3m_4m_5$,
where all prime factors $p$ of $m_i$ satisfy $\ell (p)=i$.
\end{itemize}
If $G(x)$ is the counting function of the set ${\cal G}$ then,
when $x$ tends to infinity,
\begin{equation}\label{411}
G(x) =\frac {Cx}{(\log x)^{1/4}}\left(1+{\cal O}\left(\frac {\log \log
x}{\log x}\right)\right),
\end{equation}
where
\begin{equation}\label{412}
C=\frac {H_{1,1,\xi_0}(1)C_{1,1}}{\Gamma
\left(f_{1,1}(1)\right)}=0.61568378...,
\end{equation}
$H_{1,1,\xi_0}(1)$, $C_{1,1}$ and $f_{1,1}(1)$ are defined by
(\ref {47}),(\ref {410}) and \eqref{46}.}

\medskip

{\bf Proof of Lemma 7.} We apply Lemma 6 with $y=z=1$, $\xi
=\xi_0$ and remark that $G(x)= S_{1,1,\xi_0 }(x)$. By observing
that $( 1+\frac {1}{p-1})(1-\frac {1}{p})=1$, we have
$$H_{1,1,\xi_0}(1)=\prod_{p\in {\cal P}_3}\left( 1+\frac
{1}{2(p-1)}\right) \left( 1-\frac {1}{p}\right)^{\frac
{1}{2}}=\prod_{p\in {\cal P}_3}\left( 1-\frac {1}{2p}\right)
\left( 1-\frac {1}{p}\right)^{-\frac {1}{2}}$$
$$\asymp 1.000479390466,$$ 
$$C_{1,1}=\lim_{x\rightarrow \infty }\prod_{p\in
{\cal P}_1\cup {\cal P}_2\cup {\cal P}_4\cup {\cal P}_5,~p\leq x
}\left( 1-\frac {1}{p}\right)^{-1}\prod_{p\in {\cal P}_3,~p\leq
x}\left(1-\frac {1}{p}\right)^{\frac {-1}{2}}\prod_{p\leq x}
\left(1-\frac {1}{p}\right)^{\frac {3}{4}}$$
$$\asymp 0.75410767606.$$
The numerical value of the above Eulerian products has been computed 
by the classical method already used and described in \cite{Ben5}.
Since $\Gamma \left(f_{1,1}(1)\right)=\Gamma (\frac
{3}{4})=1.225416702465...$, we get \eqref{412}. $\Box$

\medskip
{\bf Lemma 8.}  {\it We keep the notation introduced in Lemmas $6$ and $7$.
If $(y,z)\in \{(1,1),(-1,-1)\}$, we have
\begin{equation}\label{yz1}
S_{y,z,\xi_0}(x) =\frac {C\,x}{(\log x)^{1/4}}\left(1+{\cal O}
\left(\frac {\log \log x}{\log x}\right)\right),
\end{equation}
while, if $(y,z,\xi)\notin \{(1,1,\xi_0),(-1,-1,\xi_0)\}$, we have
\begin{equation}\label{Syz}
S_{y,z,\xi}(x) ={\cal O}_r\left(\frac {x}{(\log x)^{1/4+2^{-2u-3}}}\right).
\end{equation}   }

{\bf Proof.} For $y=z=1$, Formula \eqref{yz1} follows from Lemma $7$.
For $y=z=-1$ (which does not occur for $u=0$), 
it follows from \eqref{45} and by observing that the values of 
$g_{j,y,z}(s), f_{y,z}(s), H_{y,z,\xi}(s), C_{y,z}$ do not change when replacing
$y$ by $-y$ and $z$ by $-z$.

Let us define
$$M_{y,z}=\Re(f_{y,z}(1))=\frac 16 \Re(z^6(z^2+z^{-2}+\frac 12+y(z+z^{-1}))).$$
When $\xi\ne \xi_0$, \eqref{43} implies \eqref{Syz} while, if $\xi=\xi_0$,
it follows from \eqref{45} and from the inequality to be proved
\begin{equation}\label{Myz}
M_{y,z} \le \frac 34 -\frac{1}{2^{2u+3}}, \qquad (y,z)\notin \{(1,1),(-1,-1)\}.
\end{equation}
To show \eqref{Myz}, let us first recall that $z$ is a twelfth root of unity.

If $z\ne \pm 1$, $6f_{y,z}(1)$ is equal to one of the numbers
$-3/2\pm y\sqrt 3$, $-1/2\pm y$, $3/2$ so that
$$M_{y,z}\le |f_{y,z}(1)|\le \frac 16\left(\frac 32+\sqrt 3\right)
< 0.55\le \frac 34-\frac{1}{2^{2u+3}}$$
for all $u\ge 0$, which proves \eqref{Myz}.

If $z=1$ and $y\ne 1$ (which implies $u\ge 1$), we have
$$\Re y\le \cos\frac{2\pi}{2^u}=1-2\sin^2\frac{\pi}{2^u}
\le 1-2\left(\frac 2\pi\frac{\pi}{2^u}\right)^2=1-\frac{8}{2^{2u}},$$
and
$$M_{y,1}=\frac{5}{12}+\frac 13\Re y\le \frac 34-\frac{8}{3\cdot 2^{2u}} <
\frac 34-\frac{1}{2^{2u+3}}\cdot$$

If $z=-1$ and $y\ne -1$, \eqref{Myz} follows from the preceding case 
by observing that $f_{y,z}(1)=f_{-y,-z}(1)$, which completes the proof 
of \eqref{Myz}. $\Box$

\medskip

{\bf Lemma 9.} {\it Let $\cg$ be the set defined in Lemma $7$,
$\om_j$ and $\alpha'$ be the functions given by (\ref{29}) and
(\ref{alphaprime}). For $0\leq j\leq 11$, $r$, $u$, $\lambda $, $t\in
\N_0$ such that $t$ is odd, we let ${\cal G}_{j,r,u,\lambda ,t}$ be
the set of integers $n=2^{\om_3(m)}m$ in ${\cal G}$ with the
following conditions:
\begin{itemize}
  \item $\alpha'(m) \equiv j \pmod{12},$
  \item $\om_2(m)+\om_4(m)\equiv \lambda \pmod {2^{u}}$,
  \item $m\equiv t \pmod{2^{r+1}}$.
\end{itemize}
If $\rho $ is the function given by (\ref{41}), 
the counting function $G_{j,r,u,\lambda ,t}(x)$ of the set
${\cal G}_{j,r,u,\lambda ,t}$ is equal to
$$ G_{j,r,u,\lambda,t}(x) 
=\sum_{\substack{2^{\om_3(m)}m\leq x,~m\equiv t \pmod{2^{r+1}}\\ 
\alpha'(m)\equiv j \pmod{12},~\om_2(m)+\om_4(m)\equiv\lambda \pmod{2^u}\\}}
\rho (m).$$ 
If $u\ge 1$ and $\la\not\equiv j\pmod{2}$, ${\cal G}_{j,r,u,\lambda ,t}$ is 
empty while, if $\la\equiv j\pmod{2}$, when $x$ tends to infinity, 
we have
$$ G_{j,r,u,\lambda,t}(x)  = \frac {C}{6\cdot 2^{r+u}}\frac
{x}{(\log x)^{\frac {1}{4} }}\left( 1+{\cal O}\left( \frac {
1}{(\log x)^{2^{-2u-3}}}\right)\right),$$ 
where $C$ is the constant given by (\ref {412}).

\ni
If $u=0$, then
$$ G_{j,r,0,0,t}(x) = \frac {C}{12\cdot 2^r}\frac
{x}{(\log x)^{\frac {1}{4} }}\left( 1+{\cal O}\left( \frac {1}
{(\log x)^{1/8}}\right)\right),$$         }

\medskip

{\bf Proof.} If $u\ge 1$, it follows from \eqref{alphaprime} that 
$\al'(m)\equiv \om_2(m)+\om_4(m) \pmod{2}$; therefore, if 
$j\not\equiv \la \pmod{2}$, then ${\cal G}_{j,r,u,\lambda ,t}$ is empty.
Let us set 
$$\zeta =e^{\frac {2i\pi }{2^u}},~~\mu=e^{\frac {2i\pi }{12}}.$$ 
By using the relations of orthogonality: 
$$\sum_{j_2=0}^{11}\mu^{j_2\alpha'
(m)}\mu^{-jj_2}=\left\{
\begin {array}{cl}12& \hbox{if}~ \alpha'\equiv j~(\bmod ~ 12)
\\ 0 &
\hbox{\text if not,}
\end{array}\right.$$
$$\sum_{j_1=0}^{2^{u}-1}\zeta^{-\lambda j_1
}\zeta^{j_1\left(\om_2(m)+\om_4(m)\right)}=\left\{
\begin {array}{cl}2^u& \hbox{if}~\om_2(m)+\om_4(m)\equiv \lambda~(\bmod ~ 2^u)
\\ 0 &
\hbox{\text if not,}
\end{array}\right.$$
$$\sum_{\xi \bmod ~2^{r+1} }\overline{\xi }(t)\xi (m)=\left\{
\begin {array}{cl}\varphi (2^{r+1})=2^r& \hbox{if}~m\equiv t~(\bmod ~ 2^{r+1})
\\ 0 &
\hbox{\text if not,}
\end{array}\right.$$
we get 
$$G_{j,r,u,\lambda ,t}(x)=\frac {1}{12\cdot 2^{r+u}}\sum_{\xi
\bmod ~2^{r+1}}\sum_{j_1=0}^{2^u-1}\sum_{j_2=0}^{11}\overline{\xi
}(t)\zeta^{-\lambda j_1}\mu^{-jj_2}S_{\zeta^{j_1},\mu^{j_2},\xi
}(x).$$ 
In the above triple sums, the main contribution comes from
$S_{1,1,\xi_0}(x)$ and $S_{-1,-1,\xi_0}(x)$, 
and the result follows from (\ref{yz1}) and (\ref{Syz}).

If $u=0$, we have
$$G_{j,r,0,0,t}(x)=\frac {1}{12\cdot 2^{r}}\sum_{\xi \bmod ~2^{r+1}} 
\sum_{j_2=0}^{11} \overline{\xi}(t)
\mu^{-jj_2}S_{1,\mu^{j_2},\xi}(x)$$
and, again, the result follows from Lemma 8. 

\medskip

{\bf Theorem 4.} {\it Let ${\cal A}={\cal A}(1+z+z^3+z^4+z^5)$ be
the set given by (\ref {13}) and $A(x)$ be its counting function.
When $x\to\inf$, we have 
$$A(x)\sim \kappa \frac {x}{(\log x)^{\frac {1}{4}}},$$ 
where $\kappa =\frac
{74}{31}C=1.469696766...$ and $C$ is the constant of Lemma 7
defined by (\ref {412}).}

\medskip

{\bf Proof of Theorem 4.} Let us define the sets $\ca_1$, $\ca_2$,
$\ca_3$ and $\ca_4$ containing the elements $n=2^km$ ($m$ odd) of
$\ca $ with the restrictions:
\begin{eqnarray}\label{setscond}
\ca_1:&\om_3(m)\not= 0 ~\text {and}~ \om_2(m)+\om_4(m)\not=
0\nonumber\\\ca_2:& \om_3(m)\not= 0 ~\text {and}~ \om_2(m)=\om_4(m)= 0
\nonumber\\ \ca_3:& \om_3(m)= 0 ~\text {and}~ \om_2(m)+\om_4(m)\not= 0
\nonumber\\ \ca_4: & \om_2(m)=\om_3(m)=\om_4(m)= 0. \nonumber
\end{eqnarray}
We have
\begin{equation}\label{A1234}
A(x)=A_1(x)+A_2(x)+A_3(x)+A_4(x).
\end{equation}
Further, for $i=2,~3,~4$, it follows from Lemma $5$ that
$A_i(x)={\cal O}\left( \frac {x}{(\log x)^{\frac {1}{3}}}\right)$
and therefore
\begin{equation}\label{A1}
A(x)=A_1(x)+{\cal O}\left( \frac {x}{(\log x)^{\frac
{1}{3}}}\right).
\end{equation}
Now, we split $\ca_1$ in two parts $\cb$ and $\widehat{\cb}$ by
putting in $\cb$ the elements $n\in \ca_1$ which are coprime with
$31$ and in $\widehat{\cb}$ the elements $n\in\ca_1$ which are
multiples of $31$. Let us recall that, from Remark $2$, no element
of $\ca$ is a multiple of $31^2$. Therefore,
\begin{equation}\label{A1B}
A_1(x)=\cb (x)+\widehat{\cb}(x)
\end{equation}
with
\begin{equation}\label{BBhat}
\cb (x)=\sum_{n=2^km\in \ca_1,~n\leq x}\rho (m),~~\widehat{\cb}
(x)=\sum_{n=2^k31m\in \ca_1,~n\leq x}\rho (m).
\end{equation}
Let us consider $\cb (x)$; the case of $\widehat{\cb}$ will be
similar. We define
\begin{equation}\label{nui}
 \nu_i =v_2(E_i)-1=\left\{
\begin {array}{crll}-1& \hbox{if} &~ i\equiv 0,~1,~3,~4 \pmod{6}
\\ 0&\hbox{if}&~ i\equiv 2 \pmod{6}\\
2&\hbox{if}&~ i\equiv 5 \pmod{6}
\end{array}\right.
\end{equation}
so that, if $\widehat{E_i}$ is the odd part of $E_i$ 
(cf. (\ref{E}) and Table 1), we have
\begin{equation}\label{Eihat}
\widehat{E_i}=2^{-1-\nu_i}E_i.
\end{equation}
In view of Theorem $2$ (i), if $i=\al'(m) \bmod 12$ then
\begin{equation}\label{419prime}
\ga(m)-\om_3(m)=\nu_i.
\end{equation}
Further, an element $n=2^km$ ($m$ odd) belonging to $\ca_1$ is said
of index $r\geq 0$ if $k=\gamma (m)+r$. For $r\geq 0$ and $0\leq
i\leq 11$,
\begin{equation}\label{Tri}
T_r^{(i)}(x)=\sum_{\substack{n=2^{\gamma (m)+r}m\in \ca_1,~n\leq
x\\\alpha' (m)\equiv i \pmod{12}}}\rho (m)
=\sum_{\substack{n=2^{\gamma (m)+r}m\in \ca_1,~2^{\om_3(m)}m\leq
2^{-r-\nu_i}x\\\alpha' (m)\equiv i \pmod{12}}}\rho (m)
\end{equation}
will count the number of elements of $\ca_1$ up to $x$ of index
$r$ and satisfying $\alpha' (m)\equiv i \pmod{12}$, so that
\begin{equation}\label{BTi}
\cb (x)=\sum_{r\geq 0}\;\;\sum_{i=0}^{11}\;T_r^{(i)}(x).
\end{equation}
Since $\gamma (m)\geq 0$, from the first equality in (\ref{Tri}),
each $n$ counted in $T_r^{(i)}(x)$ is a multiple of $2^r$, hence
the trivial upper bound
\begin{equation}\label{Triv}
\sum_{i=0}^{11}T_r^{(i)}(x)\leq \frac {x}{2^r}\cdot
\end{equation}
Since $\nu_i\geq -1$, the second equality in (\ref{Tri}) implies
\begin{equation}\label{TriG}
\sum_{i=0}^{11}T_r^{(i)}(x)\leq G(2^{1-r}x)
\end{equation}
with $G$ defined in Lemma $7$. Moreover, from Lemma $7$, there
exists an absolute constant $K $ such that, for $x\geq 3$,
\begin{equation}\label{GK}
G(x)\leq K \frac {x}{(\log x)^{\frac {1}{4}}}\cdot
\end{equation}
Now, let $R$ be a large but fixed integer; $R'$ is defined in
terms of $x$ by $2^{R'-1}\leq \sqrt{x}<2^{R'}$ and $R"=\frac {\log
x}{\log 2}$. Since $T_r^{(i)}(x)$ is a non-negative integer, (\ref{Triv}) 
implies that $T_r^{(i)}(x)=0$ for $r>R"$. If $x$ is 
large enough, $R<R'<R''$ holds. Setting
\begin{equation}\label{BR}
\cb_R(x)=\sum_{r=0}^R\;\;\sum_{i=0}^{11}\;T_r^{(i)}(x),
\end{equation}
from \eqref{BTi}, we have 
$$\cb (x)-\cb_R(x)=S'+S",$$ 
with
$$S'=\sum_{r=R+1}^{R'}\sum_{i=0}^{11}T_r^{(i)}(x),\qquad
S"=\sum_{r=R'+1}^{R"}\sum_{i=0}^{11}T_r^{(i)}(x).$$
The definition of $R'$ and (\ref{Triv}) yield 
$$S"\leq \sum_{r=R'+1}^{R"}\frac {x}{2^r}\leq
\sum_{r=R'+1}^{\inf}\frac {x}{2^r}=
\frac {x}{2^{R'}}\leq \sqrt{x},$$ 
while (\ref{TriG}), (\ref{GK}) and the definition of $R'$ give 
\begin{eqnarray*}
S' &\leq&
\sum_{r=R+1}^{R'}G\left(\frac {x}{2^{r-1}}\right)\leq \sum_{r=R+1}^{R'}\frac
{2K x}{2^r\left( \log \frac {x}{2^{R'-1}}\right)^{\frac{1}{4}}}\\
&\leq& \frac {2^{\frac {5}{4}}K x}{(\log x)^{\frac{1}{4}}}
\sum_{r=R+1}^{R'}\frac {1}{2^r}\leq \frac {3K x}{2^R(\log x)^{\frac {1}{4}}},
\end{eqnarray*} 
so that, for $x$ large enough, we have
\begin{equation}\label{BmBR}
0\leq \cb (x)-\cb_R(x)\leq \sqrt{x}+\frac {3K x}{2^R(\log x)^{\frac {1}{4}}}\cdot
\end{equation}

We now have to evaluate $T_r^{(i)}(x)$; we shall distinguish two
cases, $r=0$ and $r\geq 1$.

\medskip

{\bf Calculation of $T_0^{(i)}(x)$.}

\medskip

\ni
From (\ref{Tri}), we have
$$T_0^{(i)}(x)=\sum_{\substack{n=2^{\ga(m)}m\in \ca_1,\;n\leq x\\
\al' (m)\equiv i \pmod{12}}} \rho (m)= 
\sum_{\substack{n=2^{\ga(m)}m\in\ca,\;n\leq x,\;\om_3\not= 0,\;
\om_2+\om_4\not= 0\\
\al' (m)\equiv i \pmod{12}}} \rho (m).$$ 
From Theorem 3, we know that $2^{\ga(m)}m\in\ca$. Hence,
$$T_0^{(i)}(x)=\sum_{\substack{2^{\ga(m)}m\leq x,
\;\om_3\not= 0,~\om_2+\om_4\not= 0\\
\al' (m)\equiv i \pmod{12}}} \rho (m),$$ 
which, by use of \eqref{419prime}, gives
$$T_0^{(i)}(x)=\sum_{\substack{2^{\om_3(m)}m\leq
2^{-\nu_i}x,~\om_3\not= 0,~\om_2+\om_4\not= 0\\
\al' (m)\equiv i \pmod{12}}} \rho (m).$$ 
But, at the cost of an error term ${\cal
O}\left( \frac {x}{(\log x)^{\frac {1}{3}}}\right)$, Lemma $5$
allows us to remove the conditions $\om_3\not= 0$, $\om_2+\om_4\not= 0$,
and to get from the second part of Lemma $9$, 
$$T_0^{(i)}(x)=G_{i,0,0,0,1}\left(\frac
{x}{2^{\nu_i}}\right)+{\cal O}\left( \frac {x}{(\log x)^{\frac
{1}{3}}}\right)$$
\begin{equation}\label{TOi}
=\frac {C}{12}\frac {x}{2^{\nu_i}(\log x)^{\frac {1}{4}}}\left(
1+{\cal O}\left( \frac {1}{(\log x)^{1/12}}\right)\right).
\end{equation}

\medskip

{\bf Calculation of $T_r^{(i)}(x)$ for $r\geq 1$.}

\medskip

\ni
Under the conditions $\om_3\not= 0$ and $\om_2+\om_4\not= 0$, 
from \eqref{mSmfinal}, \eqref{228}, \eqref{alphaprime}, \eqref{Eihat} and   
\eqref{419prime}, we get
$$Z(m)=3^{\lceil  \frac{\om_2+\om_4}{2}-1\rceil }\widehat{E}_{\alpha' (m)}.$$ 
From (\ref{Tri}), it follows that
$$T_r^{(i)}(x)=\sum_{\substack{n=2^{\ga(m)+r}m\in\ca,\; n\leq x,\;
\om_3\not= 0,~\om_2+\om_4\not= 0\\
\al' (m)\equiv i \pmod{12}}} \rho (m).$$ 
Now, by Theorem 3, we know that $2^{\gamma(m)+r}m$ 
belongs to $\ca$ if there is some $l\in\cs_r=\{2^r+1, ...,
2^{r+1}-1\}$ such that $m\equiv l^{-1}Z(m) \bmod ~2^{r+1}$. Note
that the order of $3$ modulo $2^{r+1}$ is $2^{r-1}$
if $r\ge 2$ and $2^r$ if $r=1$. We choose
$$u=r+1$$
so that $\om_2+\om_4\equiv \la \pmod{2^{r+1}}$ implies 
$3^{\lceil \frac \la2-1 \rceil}\equiv 
3^{\lceil \frac{\om_2+\om_4}{2}-1 \rceil} \pmod{2^{r+1}}$.
Therefore, we have 
$$T_r^{(i)}(x)=\sum_{l\in \cs_r}\sum_{\lambda
=0}^{2^{r+1}-1}~~\sum_{\substack{2^{\om_3(m)}m\leq
2^{-\nu_i-r}x,~\om_3\not= 0,~\om_2+\om_4\not= 0\\ \alpha' (m)\equiv i
\pmod{12},~\om_2+\om_4 \equiv \lambda \pmod{2^{r+1}}\\m\equiv
l^{-1}3^{\lceil \frac {\lambda }{2}-1\rceil }\widehat{E_i} 
\pmod{2^{r+1}}}}\rho (m).$$ 
As in the case $r=0$, we can remove the
conditions $\om_3\not= 0$ and $\om_2+\om_4\not= 0$ in the last sum by
adding a ${\cal O}\left( \frac {x}{(\log x)^{\frac
{1}{3}}}\right)$ error term, and we get by Lemma $9$ for $r$ fixed
$$T_r^{(i)}(x)=\sum_{l\in \cs_r}\sum_{\substack{\lambda=0\\\la\equiv i ~(\bmod 2)}}
^{2^{r+1}-1}G_{i,r,r+1,\lambda ,l^{-1}3^{\lceil \frac {\lambda
}{2}-1\rceil }\widehat{E_i}}\left( \frac {x}{2^{\nu_i+r}}\right)
~+~{\cal O}\left( \frac {x}{(\log x)^{\frac {1}{3}}}\right)$$
\begin{equation}\label{Trix}
=\frac {C}{24}\frac {x}{2^{\nu_i+r}(\log x)^{\frac {1}{4}}}\left(
1+{\cal O}\left(\frac {1}{(\log x)^{2^{-2r-5}}}\right)\right).
\end{equation}
From (\ref{BR}), (\ref{TOi}), (\ref{Trix}) and (\ref{nui}), we
have 
\begin{eqnarray*}
\cb_R(x)&=&\frac {C x}{12(\log x)^{\frac {1}{4}}}\left(
\left(\sum_{i=0}^{11}\frac {1}{2^{\nu_i}}\right)
\left( 1+\frac{1}{2}\sum_{r=1}^R\frac {1}{2^r}\right)
+{\cal O}\left( \frac {1}{(\log x)^{2^{-2R-5}}}\right)   \right)\\
&=& \frac {37}{24}\frac
{C x}{(\log x)^{\frac {1}{4}}}\left( \frac {3}{2}-\frac
{1}{2^R}\right)\left( 1+{\cal O}\left( \frac {1}{(\log
x)^{2^{-2R-5}}}\right)\right).
\end{eqnarray*}
By making $R$ going to infinity, the above
equality together with (\ref{BmBR}) show that
\begin{equation}\label{Bx}
\cb (x)\sim \frac {37}{16}\frac {C x}{(\log x)^{\frac
{1}{4}}},~~x\rightarrow \infty.
\end{equation}
In a similar way, we can show that $\widehat{\cb}(x)$ defined in
(\ref{BBhat}) satisfies $$\widehat{\cb}(x)\sim \frac {1}{31}\cb
(x)\sim \frac {37}{16\cdot 31}\frac {x}{(\log x)^{\frac {1}{4}}}$$
which, with (\ref{A1B}) and (\ref{A1}), completes the proof of
Theorem $4$ with 
$$\kappa =\frac {37}{16}\left( 1+\frac
{1}{31}\right) C=\frac {74}{31}C=1.469696766....$$ 

\medskip

{\bf Numerical computation of $A(x)$.}

\ni
There are three ways to compute $A(x)$. The
first one uses the definition of $\ca$ and simultaneously
calculates the number of partitions $p({\cal A}, n)$ for $n\leq
x$; it is rather slow. The second one is based on the relation
(\ref{110}) and the congruences (\ref{sigmaequivmod2k+1}) and
(\ref{Sm31}) satisfied by $\sigma ({\cal A},n)$. The third one
calculates $\om_j(n),~0\leq j\leq 5$, in view of applying Theorem
$1$. The two last methods can be encoded in a sieving process

The following table displays the values of $A(x)$, $A_1(x)$, ...,
$A_4(x)$ as defined in (\ref{A1234}) and also 
$$c(x)=\frac
{A(x)(\log x)^{\frac {1}{4}}}{x},~~c_1(x)=\frac {A_1(x)(\log
x)^{\frac {1}{4}}}{x}.$$ 
It seems that $c(x)$ and $c_1(x)$  converge very slowly to
$\kappa=1.469696766\ldots$, which is impossible to guess from the table.

\begin{center}
$$
\begin{array}{|l|l|l|l|l|l|l|l|}
  \hline
  x & A(x) & c(x) & A_1(x) & c_1(x) & A_2(x) & A_3(x) & A_4(x)  \\
  \hline
  10^3 & 480 & 0.7782 & 20 & 0.032 & 44 & 233 & 183 \\
  10^4 & 4543 & 0.7914 & 361 & 0.063 & 532 & 2294 & 1356 \\
  10^5 & 43023 & 0.7925 & 5087 & 0.094 & 5361 & 21810 & 10765 \\
  10^6 & 411764 & 0.7939 & 60565 & 0.117 & 52344 & 208633 & 90222 \\
  10^7 & 3981774 & 0.7978 & 680728 & 0.136 & 506199 & 2007168 & 787679 \\
  10^8 & 38719773 & 0.8022 & 7403138 & 0.153 & 4887357 & 19390529 & 7038749 \\
\hline
\end{array}
$$
\end{center}

\medskip

{\bf Thanks}

\ni
We are pleased to thank A. \sark who first considered the sets $\ca$'s 
such that the number of partitions $p(\ca,n)$ is even for $n$ large enough
for his interest in our work and X. Roblot for valuable discussions 
about $2$-adic numbers.

\def\refname{References}

\end{document}